 \documentclass[12pt]{article} \usepackage{amsmath,theorem,amssymb,amscd}

\usepackage{rotating,graphicx}

  \newtheorem{lemma}[equation]{Lemma}  \newtheorem{theorem}[equation]{Theorem} \theoremstyle{plain}
\theorembodyfont{\rmfamily}  \newtheorem{definition}[equation]{Definition}   \newenvironment{pf}{\noindent \textbf{Proof.}}{\hfill{$\square$}\\} \numberwithin{equation}{section}

\begin{document}

\date{}
\title{Terminal Resolutions of Brauer Pairs}
\author{Basil Nanayakkara\\basil.n@unb.ca}

\maketitle
\thispagestyle{empty}

\begin{abstract}
A \emph{Brauer pair} is a pair $(X,\alpha)$ where $X$ is a quasi-projective variety over an algebraically closed field and $\alpha$ is an element in the 2-torsion part of the Brauer group of the function field of $X$.  A Brauer pair $(Y,\alpha )$ is a \emph{terminal pair} if the Brauer discrepancy of $(Y,\alpha)$ is positive.  We show that given a Brauer pair $(X,\alpha)$, there is a terminal pair $(Y,\alpha )$ with a birational morphism $Y \longrightarrow X.$  In short, any Brauer pair admits a terminal resolution.
\end{abstract}
\maketitle

\begin{section}{Introduction}
The problem considered in this article is central to research in the area of maximal orders on algebraic varieties.   The works of Daniel Chan, Colin Ingalls and Rajesh Kulkarni (\cite{CIH}, \cite{CK1}, \cite{CK2}, \cite{I}), considered various questions regarding orders on algebraic surfaces.  The main idea was to build a minimal model program for maximal orders on surfaces.  This was  successfully established in a seminal article in the area by Daniel Chan and Colin Ingalls \cite{CI}.  In this article, they define the notion of terminal orders (which are analogs of smooth surfaces) and then prove the main theorem:  any maximal order has a terminal resolution.  This led to several articles in which terminal models of orders on surfaces were classified.  This classification has been a significant achievement of the past decade.

In Section 2, we explain how an element in the Brauer group of the function field of a variety induces a boundary divisor, via a complex that appears in the coniveau spectral sequence of the variety.  In Section 3, we describe birational geometry of Brauer pairs after introducing the notion of Brauer discrepancy of a pair.  Finally, in Section 4, we prove the main result that any Brauer pair admits a terminal resolution.


\end{section}

\section{ Logarithmic pairs from Brauer pairs }
\label{conispec}
Level 1 terms of the coniveau spectral sequence for a smooth algebraic variety $X$ over an algebraically closed field $k$ were written by Grothendieck \cite{Groth}:  One has 
$$
E_1^{i,j}= \bigoplus _{x \in X^{(i)}} H^{j-i}(k(x),\mu_n^{\otimes(1-i)})
$$ 
where $\mu_n$ is the group of $n^{\mathrm{th}}$ roots of unity in $k$, and $X^{(i)}$ is the set of all irreducible subvarieties of $X$ of codimension $i$.  The cohomology mentioned is Galois cohomology.  The tensor product is over ${\bf Z}/{n {\bf Z}}.$   By definition, $\mu_n^{-1}:=\mathrm{Hom}(\mu_n, {\bf Z}/{n \bf Z}),$ and we write $\mu_n^{\otimes(-m)}$ for $(\mu_n^{-1})^{\otimes m}$ when $m$ is positive. For more details, see Section 3.5 of \cite{LeB}.

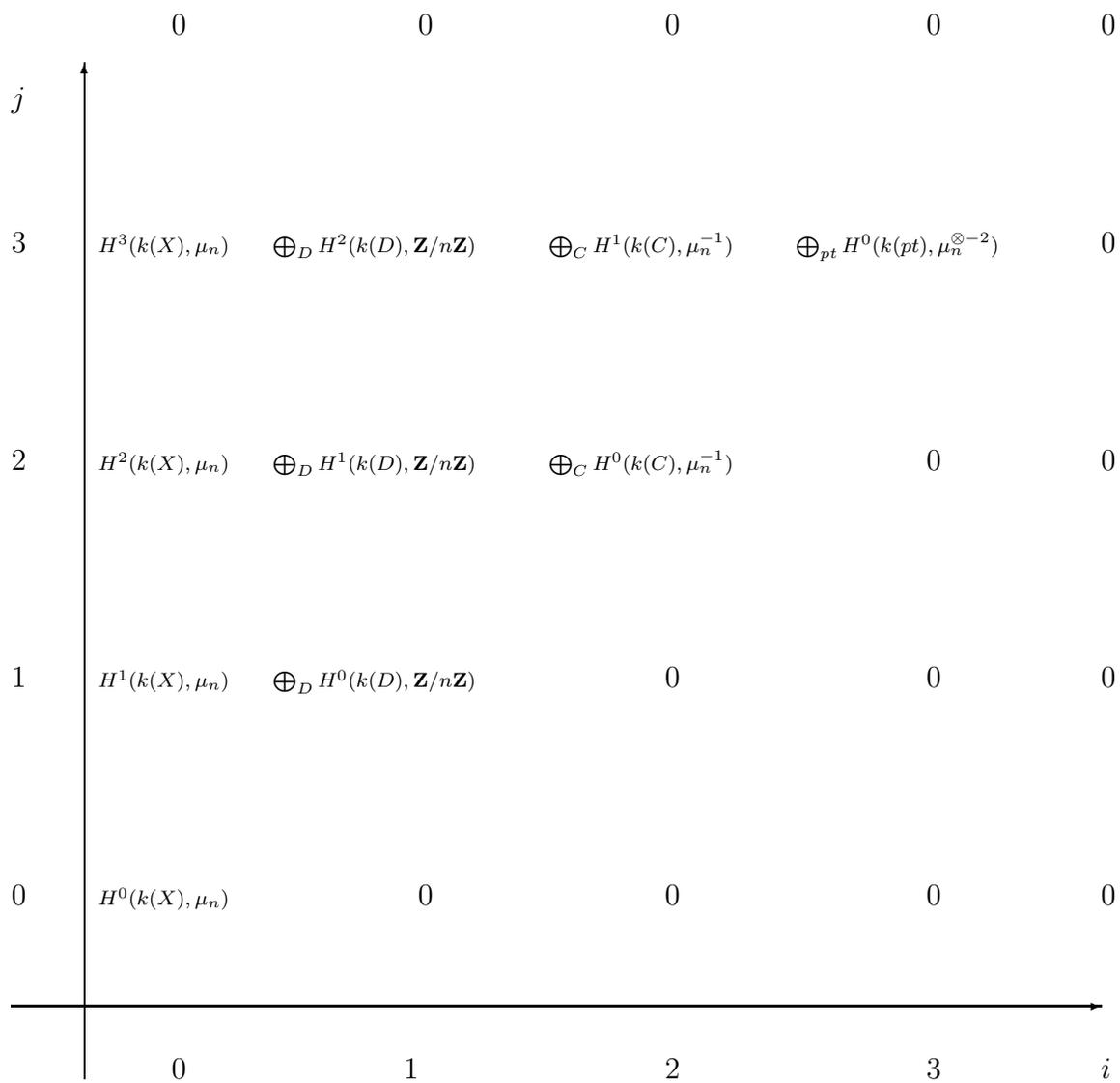
\begin{figure}[t]
\setlength{\unitlength}{2cm}
\centering
\caption{Level 1 of the coniveau spectral sequence for an irreducible 3-fold $X$}
\begin{picture}(5,10)(1,-1)

\put(-7,71){$\underline{C1}$}
\put(0,0.5){\vector(0,1){7}}

\put(.6,.5){$0$}
\put(2.2,.5){$1$}
\put(4,.5){$2$}
\put(5.8,.5){$3$}
\put(7,.5){$i$}

\put(-.5,1.7){$0$}
\put(-.5,3.2){$1$}
\put(-.5,4.7){$2$}
\put(-.5,6.2){$3$}
\put(-.5,7.2){$j$}

\put(-.5,1){\vector(1,0){7.5}}

\put(.6,7.7){$0$}
\put(.1,1.7){\scriptsize $H^0(k(X),\mu_n)$}
\put(.1,3.2){\scriptsize $H^1(k(X),\mu_n)$}
\put(.1,4.7){\scriptsize $H^2(k(X),\mu_n)$}
\put(.1,6.2){\scriptsize $H^3(k(X),\mu_n)$}

\put(2.3,1.7){$0$}
\put(1.3,3.2){\scriptsize $\bigoplus_D H^0(k(D),{\bf Z}/{n \bf Z})$}
\put(1.3,4.7){\scriptsize $\bigoplus_D H^1(k(D),{\bf Z}/{n \bf Z})$}
\put(1.3,6.2){\scriptsize $\bigoplus_D H^2(k(D),{\bf Z}/{n \bf Z})$}
\put(2.3,7.7){$0$}

\put(4,1.7){$0$}
\put(4,3.2){$0$}
\put(3.2,4.7){\scriptsize $\bigoplus_C H^0(k(C),\mu_n^{-1})$}
\put(3.2,6.2){\scriptsize $\bigoplus_C H^1(k(C),\mu_n^{-1})$}
\put(4,7.7){$0$}

\put(5.8,1.7){$0$}
\put(5.8,3.2){$0$}
\put(5.8,4.7){$0$}
\put(4.9,6.2){\scriptsize $\bigoplus_{pt} H^0(k(pt),\mu_n^{\otimes -2})$}
\put(5.8,7.7){$0$}

\put(7,1.7){$0$}
\put(7,3.2){$0$}
\put(7,4.7){$0$}
\put(7,6.2){$0$}
\put(7,7.7){$0$}

\put(3.1,61){$\oplus$}
\put(-1.2,66){$\oplus$}
\put(-6,59){x}
\put(6,59){y}

\end{picture}

\label{coniveau}  
\end{figure}

Accordingly, if $X$ is an irreducible 3-fold, we get Figure \ref{coniveau} on page \pageref{coniveau} for the first quadrant of level 1.  In that figure,  D, C, pt are prime divisors, irreducible curves and points of $X$ respectively.  The row j=2 has the same form for any irreducible n-fold where $n \ge 2$ with the interpretation that $C$ and pt are irreducible subvarieties of codimension 2 and 3 respectively.

Now, from row $j=2$, we obtain the complex
$$
H^2(k(X), \mu_n) \longrightarrow \bigoplus _{D} H^1(k(D), {\bf Z}/{n \bf Z}) \longrightarrow \bigoplus _{C} H^0(k(C), \mu_n^{-1}) \longrightarrow 0.
$$

We know that  $H^2(k(X), \mu_n) \cong Br_n(k(X))$, where $ Br_n(k(X))$ is the n-torsion part of the Brauer group $Br (k(X))$ of $k(X)$. (See section 4.4 of \cite{GS}).  Therefore, we get the complex
$$
Br_n(k(X)) \longrightarrow \bigoplus _{D} H^1(k(D), {\bf Z}/{n \bf Z}) \longrightarrow \bigoplus _{C} \mu_n^{-1} \longrightarrow 0.
$$

This tells us, in particular, that any element $\alpha$ in $Br_2(k(X))$, induces a (possibly ramified) 2-sheeted cover or a 1-sheeted cover on each irreducible divisor $D$. Note that the ramifications must cancel on the irreducible subvarieties $C$, since the sequence above is a complex.  This simple observation will play an important role in our study.



\begin{definition}
A \emph{Brauer pair} is a pair $(X, \alpha)$ where $X$ is a quasiprojective variety over an algebraically closed field and $\alpha$ is an element in $Br_2  (k(X))$.
\end{definition}

Let $(X,\alpha)$ be a Brauer pair.   Then $\alpha$ induces a boundary divisor $\Delta_{X,\alpha}$ on $X$ as follows:

Consider the complex
\begin{equation}
\label{complex}
\mathrm{Br}_n \ k(X) \overset{a}{\longrightarrow} \bigoplus _{D} H^1(k(D),{\bf Z}/n {\bf Z}) \longrightarrow \bigoplus _{C} \mu_n^{-1} \longrightarrow 0
\end{equation} 
that we obtained above, where $D$ runs through all the irreducible divisors of $X$ and $C$ runs through all the irreducible subvarieties of codimension 2 of $X$.  For a given irreducible divisor $D$, let $a(\alpha)_D$ be the image of $\alpha$ indexed by $D$.  Since $H^1(k(D),{\bf Q/Z})$ classifies cyclic covers of $D$,  $a(\alpha )_D$ determines a ramified cover of $D$.  Let $e_D$ be the degree of this cover.  We define the boundary divisor $\Delta _{X, \alpha}$ to be,
$$
\Delta _{X, \alpha} := \sum _D \left(1-\frac{1}{e_D}\right)D
$$ 
where $D$ runs through all the prime divisors of $X$ such that $a(\alpha)_D \ne 0$. (See Section 3.3 of \cite{CI}).

\begin{section}{Birational Geometry of Brauer Pairs}
\label{BGofBP}


In the following, by a \textit{divisor over $X$} we mean an irreducible divisor $E \subseteq Y$ where $Y$ is a normal variety with a birational morphism $Y \longrightarrow X$.

Two divisors $D_1, D_2$ of $X$ are said to be \emph{numerically equivalent} if $D_1 \cdot C = D_2 \cdot C$ for all irreducible curves $C \subseteq X.$

In the definition below, by $K_{Z,\alpha}$ we  mean $K_Z+\Delta_{Z,\alpha}$, and $\equiv$ denotes numerical equivalence.

\begin{definition}
Let $E$ be an irreducible exceptional divisor over $X$ and $\alpha \in Br \ k(X)$.  The \emph{Brauer discrepancy} of the pair $(X, \alpha)$ along $E$, denoted by $b(E,X,\alpha)$, is the coefficient of $E$ in the formula
$$
K_{Y,\alpha} \equiv f^*K_{X,\alpha}+\sum_i b(E_i,X,\alpha)E_i
$$
where $E$ is one of the exceptional divisors $E_i$ of $f$ where $f:Y \longrightarrow X$ is a birational morphism.  
\end{definition}

\begin{definition}
Let $f:Y \longrightarrow X$ and $g:Y' \longrightarrow X$ be birational morphisms.   Suppose $E \subseteq Y$ and $E' \subseteq Y'$ are $f$-exceptional and $g$-exceptional divisors respectively.   We say that $E$ and $E'$ are \emph{isomorphic as divisors} if there exist open sets $U, U'$  containing $E$ and $E'$ respectively and an isomorphism $\phi:U \longrightarrow U'$ such that $\phi|_E:E \longrightarrow E'$ is an isomorphism and $\phi=g^{-1}\circ f$ on $U \setminus \{f\mbox{-exceptional divisors} \}.$
\end{definition}

The following lemma is the reason why we suppress the birational morphism $f$ and the variety $Y$ in the notation $b(E,X, \alpha )$ for the Brauer discrepancy of the pair $(X,\alpha )$ along $E.$  The analogous remark for usual discrepancy (see below for the definition), is mentioned in Remark 2.23 of \cite{KM}.

\begin{lemma}
$b(E,X,\alpha)$ does not depend on the particular birational morphism.
\end{lemma}
\begin{pf}
Suppose $f:Y \longrightarrow X$ and $g:Y' \longrightarrow X$ are birational morphisms, $E \subseteq Y$ and $E' \subseteq Y'$ are $f-$exceptional and $g-$exceptional divisors respectively, that are isomorphic as divisors.  Then there is an isomorphism  $\phi:U \longrightarrow U'$  on a neighborhood $U$ of $E$ such that $\phi=g \circ f^{-1}$ almost everywhere on $U$, and $U'$ is a neighborhood of $E'.$  Now, we have,
$$
K_{Y,\alpha} \equiv f^*K_{X,\alpha}+\sum_i b(E_i,X,\alpha)E_i +b(E,X,\alpha)E
$$
where $E_i$ are $f-$exceptional divisors, $E_i \ne E$.

Similarly,
$$
K_{Y',\alpha} \equiv g^*K_{X,\alpha}+\sum_j b(E'_j,X,\alpha)E'_j +b(E',X,\alpha)E'
$$
where $E'_j$ are $g-$exceptional divisors, $E'_j \ne E'$.

Now,
$$
(K_{Y,\alpha}-f^*(K_{X,\alpha}))\vert _U = \phi ^*((K_{Y',\alpha} -g^*(K_X,\alpha))\vert _{U'}).
$$
Therefore,
$$
\sum b(E_i, X, \alpha)E_i\vert _U +b(E,X,\alpha )E = \phi ^* \left(\sum b(E'_j, X, \alpha)E'_j\vert _U\right) +b(E', X, \alpha )E.
$$
Hence $b(E,X,\alpha)=b(E',X,\alpha).$

\end{pf}

For comparison purposes, and also for later use, we give the definition of the (usual) discrepancy of a log pair $(X,\Delta)$, which chronologically preceded the notion of Brauer discrepancy.

\begin{definition}
Let $(X,\Delta)$ be a logarithmic pair and $E$ an exceptional divisor over $X.$  The \emph{discrepancy} of the pair $(X,\Delta)$ along  $E$, denoted by $a(E,X,\Delta)$, is the coefficient of $E$ in
$$
K_Y+f_*^{-1}\Delta \equiv f^*(K_X+\Delta)+\sum_{E_i}a(E_i,X,\Delta)E_i
$$
where $E$ is one of the exceptional divisors $E_i$ and $f:Y \longrightarrow X$ is a birational morphism.
\end{definition}

This is the usual notion of discrepancy in algebraic geometry.

The lemma below tells us how Brauer discrepancy and the (usual) discrepancy are related.

\begin{lemma}
Let $(X,\alpha)$ be a Brauer pair, $E$ an exceptional divisor over $X$, and $\Delta_{X,\alpha}$ the boundary divisor induced by $\alpha$.  Then,
$$
b(E,X,\alpha)=a(E,X,\Delta_{X,\alpha})+1-\frac{1}{e}
$$
where $e$ is the degree of the cover on $E$ induced by $\alpha$.
\end{lemma}
For a proof, see the proof of Proposition 3.15 of \cite{CI}.

We define the \emph{Brauer discrepancy} of the Brauer pair $(X,\alpha)$ to be,
$$
\mbox{bdiscrep}(X,\alpha) := \inf \{e(E,X,\alpha) \cdot b(E,X,\alpha): E \mbox{ is an exceptional divisor over } X\},
$$
where $e(E,X,\alpha)$ is the degree of the cover on $E$ induced by $\alpha$.

\begin{definition}
A pair $(X,\alpha)$ is a \emph{terminal pair} if bdiscrep$(X,\alpha) > 0$.
\end{definition}

\begin{definition}
A \emph{ terminal resolution } of $(X,\alpha)$ is a birational morphism $Y \longrightarrow X$ from a smooth variety $Y$ such that the pair $(Y,\alpha)$ is a terminal pair.
\end{definition}

Now we can state our main theorem.
\begin{theorem}
Any Brauer pair $(X,\alpha)$  has a terminal resolution.
\end{theorem}


We will see that we can arrive at a terminal pair by successively blowing up $(X, \Delta_{X, \alpha})$. 


Using Hironaka's theorem on resolution of singularities \cite{Hiro}, we can improve $(X,\Delta_{X,\alpha})$ such that $X$ is smooth and $\Delta_{X,\alpha}$ is simple normal crossing.  Thus, \'etale-locally, the improved $X$ has the form $X=\mbox{spec } k\{x_1, x_2, \dots , x_n\}.$  The boundary divisor is of the form,
$$
\Delta_{X,\alpha}=\sum_{i=1}^n \left(1-\frac{1}{e_i}\right) V(x_i)
$$
where $V(x_i)$ is the hyperplane defined by $x_i=0.$  \hfill \\
I. e. $V(x_i)= \mbox{spec }k\{x_1, \dots , x_{i-1}, x_{i+1}, \dots , x_n \} $  with $V(x_i) \hookrightarrow X$ the dual of the map $k\{x_1, \dots , x_n\} \rightarrow k\{x_1, \dots , x_{i-1}, x_{i+1}, \dots , x_n\}$ arrived at by setting $x_i=0.$  The number $e_i$ is the degree of the cover on $V(x_i)$ induced by $\alpha.$  Since $\alpha$ is in the 2-torsion part of $Br \ k(X)$, we have $e_i \in \{1,2\}$, by complex \ref{complex} of Section \ref{conispec}.

 Suppose $D_1$ and $D_2$ are prime divisors of $X$ with covers on them.  Let $C=D_1 \cap D_2$.  If one of the covers ramifies on $C$ then the other also must ramify on $C$, since the sequence (2.2) is a complex.

\section{Proof of the main result}
In this section we will prove the main result that given a Brauer pair $(X,\alpha)$, there is a smooth variety $Y$ and a birational morphism $Y\longrightarrow X$ such that $(Y,\alpha)$ is a terminal pair.


\vskip 5mm
First, we present a few lemmas.

\begin{lemma} 
\label{bdiscrep1}
Let $X$ be a smooth variety, $Z \subseteq X$ a smooth subvariety of codimension $c$ and $\alpha \in \mbox{Br }k(X).$  Suppose $p:B_Z X=Y \longrightarrow X$ is the blow-up of $X$ along $Z$, $E$ is the exceptional divisor.  Then,
$$
b(E,X,\alpha)=c-\frac{1}{e}-\sum_i a_i \cdot \mathrm{mult}_Z D_i
$$
where $e$ is the degree of the cover on $E$ determined by $\alpha$ and $\sum_i a_i D_i$ is the boundary divisor on $X$ determined by $\alpha$. 
\end{lemma}

\begin{pf}
By the definition of Brauer discrepancy,
\begin{eqnarray*}
b(E,X,\alpha)E &= &K_Y +\Delta_{Y,\alpha}-p^*(K_X+\Delta_{X,\alpha}) \\
               &=&K_Y-p^*K_X+\Delta_{Y,\alpha}-p^*\Delta_{X,\alpha}.
\end{eqnarray*}
Now, $K_Y-p^*K_X=(c-1)E$ and 
\begin{eqnarray*}
\Delta_{Y,\alpha}-p^*(\Delta_{X,\alpha}) &=& \left(1-\frac{1}{e}\right)E +p^{-1}_*(\Delta_{X,\alpha})-p^*\Delta_{X,\alpha} \\
&=&\left(1-\frac{1}{e}\right)E+p_*^{-1}\left(\sum_i a_iD_i\right)-\sum_i a_i (p^*D_i) \\
&=&\left(1-\frac{1}{e}\right)E -\sum_i a_i(p^*D_i-p_*^{-1}D_i) \\
&=&\left(1-\frac{1}{e}\right)E-\sum_i a_i \cdot (\mathrm{mult}_ZD_i)E.
\end{eqnarray*}
Thus,
\begin{eqnarray*}
&b(E,X,\alpha)E &= (c-1)E+\left(1-\frac{1}{e}\right)E-\left(\sum_ia_i \cdot \mathrm{mult}_ZD_i \right)E, \quad \mbox{ and so } \\
 &b(E,X,\alpha) &= c-\frac{1}{e}-\sum_i a_i \cdot \mathrm{mult}_ZD_i.
\end{eqnarray*}
\end{pf}

\vskip 5mm

The following lemma, which is the analog of Lemma \ref{bdiscrep1}   for usual discrepancy, appears as Lemma 2.29 in \cite{KM} without proof.

\begin{lemma}[Analogous to Lemma \ref{bdiscrep1}]
\label{discrep1}
Let $X$ be a smooth variety and $\sum a_iD_i$ is a boundary divisor on $X$.  Let $Z \subseteq X$ be a smooth subvariety of codimension $c$.  Suppose $p:B_ZX=Y \longrightarrow X$ is the blow up of $X$ along $Z$, and $E$ denotes the exceptional divisor.  Then,
$$
a(E,X,\Delta)=c-1-\sum_ia_i\cdot mult_ZD_i.
$$
\end{lemma}
Proof is similar to the proof of Lemma \ref{bdiscrep1}, and therefore omitted.

\vskip 5mm

Lemma 2.45 of \cite{KM} tells us that any exceptional divisor over a variety $X$ can be reached by finitely many blow ups.  This encourages the following definitions.
\begin{definition}
Let $E$ be an exceptional divisor over a variety $X$.  The divisor $E$ is called a \emph{level $n$ exceptional divisor}, if 
$$
n=\mbox{inf}\left\{m \in {\bf Z}^+ : E \mbox{ can be reached by $m$ successive blow ups starting from } X\right\}.
$$
\end{definition}

\begin{definition}
 A Brauer pair $(X,\alpha)$ is called \emph{level n Brauer terminal} if $b(E,X,\alpha) > 0$ for all level $m$ exceptional divisors $E$ over $X$ for $1 \le m \le n.$
\end{definition}

\begin{subsection}{Blowing up along a subvariety of codimension $\ge 3$}

In this section we show that if $E$ is an exceptional divisor generated by a blow up of a subvariety of codimension $\ge 3$, then $b(E,X, \alpha)$ is positive.

\begin{lemma}
Let $(X,\alpha)$ be an n-dimensional Brauer pair and $Z \subseteq X$ be a subvariety of codimension $c$.  Let $E$ be the exceptional divisor generated in the blow up of $X$ along $Z$.  If $c \ge 3$, then $b(E,X,\alpha)$ is positive.
\end{lemma}
\begin{pf}
 As described towards the end of Section \ref{BGofBP}, \'etale locally $(X,\Delta_{X,\alpha})$ is of the form $\left(\mbox{spec } k\{x_1, \dots , x_n\}, \sum_i \left(1-\frac{1}{e_i}\right)V(x_i)\right).$   Thus, by Lemma \ref{bdiscrep1},
$$
b(E,X,\alpha)=c-\frac{1}{e_E}-\sum_i \left( 1-\frac{1}{e_i}\right)\cdot mult_ZV(x_i).
$$

Since $Z$ is of codimension $c$, it lies on at most $c$ of $V(x_i)$.  Thus,
$$
\sum_i \left(1-\frac{1}{e_i}\right) \cdot mult_ZV(x_i) \le c-\left( \frac{1}{e_{i_1}} + \frac{1}{e_{i_2}} + \cdots + \frac{1}{e_{i_c}}\right).
$$

Therefore,
\begin{eqnarray*}
b(E,X,\alpha) & \ge & c-\frac{1}{e_E}-\left\{c-\left(\frac{1}{e_{i_1}} + \cdots +\frac{1}{e_{i_c}} \right) \right\} \\
		& \ge & \frac{1}{e_{i_1}} + \cdots +\frac{1}{e_{i_c}} - \frac{1}{e_E}\\
		& \ge & c \cdot \frac{1}{2} - \frac{1}{e_E} \quad \quad \mbox{ since } e_{i_j} =1 \mbox{ or } 2 \\
		& > & 0 \quad \mbox{ since } e_E=1 \mbox{ or } 2 \mbox{ and } c \ge 3.
\end{eqnarray*}

\end{pf}

\end{subsection}

\begin{subsection}{Blowing up along a subvariety of codimension 2}

Here we show that if $X$ is blown up along a subvariety of codimension 2, then either the Brauer discrepancy of $(X,\alpha )$ along the resulting exceptional divisor is positive, or $X$ can be blown up to obtain $X'$ so that $(X',\alpha)$ is level 1 Brauer terminal.

\begin{lemma}
\label{levelone}
 Let $(X,\alpha)$ be an n-dimensional Brauer pair with $\alpha ^2 =1$ and $Z \subseteq X$ be a subvariety of codimension $2.$  Let $E$ be the $f$-exceptional divisor where $f:Bl_Z(X) \longrightarrow X$ is the blow up of $X$ along $Z$.   Then, either $(1) \quad b(E,X,\alpha) > 0$, or $(2) \quad X$ can be blown up to obtain a level 1 Brauer terminal pair $(X',\alpha)$.
\end{lemma}

\begin{pf}
 
Locally, we have $\left( \mbox{spec } k\{x_1, \dots, x_n \}, \sum_i \left( 1-\frac{1}{e_i}\right) V(x_i) \right) $.  Since $Z$ is of codimension 2, it lies on at most two of $V(x_i)$.  If $Z$ lies on none or exactly one of $V(x_i)$, then by Lemma \ref{bdiscrep1},
$$
b(E,X,\alpha) = 2 - \frac{1}{e_E}
$$
or
$$
b(E,X,\alpha ) =2-\frac{1}{e_E} -\left( 1-\frac{1}{e_i}\right)= 1-\frac{1}{e_E} + \frac{1}{e_i}
$$
respectively.  In either case, we have $b(E,X,\alpha) >0$, because $e_E, e_i \in \{1,2\}$.

Now, suppose $Z$ lies on two of $V(x_i)$.  Assume, without loss of generality, that $Z$ lies on $V(x_1)$ and $V(x_2)$.   Then, again by Lemma \ref{bdiscrep1},
\begin{eqnarray*}
b(E,X,\alpha) & = & 2 -\frac{1}{e_E} -\left( 1-\frac{1}{e_1}+1-\frac{1}{e_2} \right) \\
		& = & \frac{1}{e_1} + \frac{1}{e_2} - \frac{1}{e_E}.
\end{eqnarray*}
Thus, we see that $b(E,X,\alpha)$ is positive except when $e_1=e_2=2$ and $e_E=1$, in which case  $b(E,X,\alpha) = 0.$   But this latter situation can occur only when the covers on $V(x_1)$ and $V(x_2)$ are double covers that do not ramify on $V(x_1,x_2)$.  (If the double covers ramify, then the induced cover on $E$ must be double cover which makes $e_E=2$, a contradiction.)  Whenever this situation occurs, we can obtain a new variety $X'$ by blowing up along $V(x_1,x_2)$.  The resulting variety $X'$ will not have two prime divisors with double covers on them without the covers ramifying on the intersection of the prime divisors.  Therefore, in the new variety this undesirable situation $(e_1=e_2=2$ and $e_E=1)$ does not occur.  Hence $(X',\alpha)$ is level 1 Brauer terminal.

\end{pf}

\end{subsection}

\begin{subsection}{Completion of the Proof}
 
In this section we show that if $(X,\alpha)$ is level 1 Brauer terminal, then it is indeed Brauer terminal and complete the proof that any Brauer pair admits a terminal resolution.

In the following lemma, $a(E,X, \Delta)$ denotes the (usual) discrepancy of the logarithmic pair $(X, \Delta)$ along an exceptional divisor $E$ over $X$.

\begin{lemma}
\label{usualdiscrepancy}
Let $(X,\alpha)$ be a Brauer pair such that $X$ is smooth and $\Delta_{X,\alpha}$ is simple normal crossing.  Let $Z \subseteq X$ be an irreducible subvariety of codimension $c$, where $c \ge 2.$   Suppose $p:B_ZX \longrightarrow X$ is the blow up of $X$ along $Z$ and $E \subseteq B_ZX$ the exceptional divisor.  Then, $a(E,X,\Delta_{X,\alpha}) \ge 0.$
 
\end{lemma}

\begin{pf}
 \'Etale locally, we have
$$
\Delta_{X,\alpha} = \sum_{i=1}^n \left(1-\frac{1}{e_i}\right) V(x_i).
$$
Applying Lemma \ref{discrep1}, we see that
$$
a\left(E,X,\Delta_{X,\alpha}\right) = c-1 -\sum_{i=1}^n \left(1-\frac{1}{e_i}\right) \cdot mult_Z V(x_i).
$$
Since $Z$ is of codimension $c$, it lies on maximum of $c$ prime divisors $V(x_i)$.  Thus,
$$
a(E,X,\Delta_{X,\alpha}) \ge c-1 -c \left( 1-\frac{1}{e_i}\right) =\frac{c}{e_i} -1 \ge 0, 
$$
since $c \ge 2$ and $e_i \in \{1,2\}$.
\end{pf}

The following lemma can be considered as a Brauer version of the composition of Lemma 2.29 and Lemma 2.30 of \cite{KM}.

\begin{lemma}
\label{composition}
Let $f:Y\longrightarrow X$ be a birational morphism, $E$ is an $f$-exceptional divisor in $Y$, $E_0$ an irreducible subvariety of $E$.  Suppose $g:Z \longrightarrow Y$ is the blow up of  $Y$ along $E_0$ and $F$ the $g$-exceptional divisor.  Suppose $b(E',X,\alpha) \ge 0$ for all $f$-exceptional divisors $E'$ and $a(F,Y,\Delta_{Y,\alpha}) \ge  0.$  Let $\lambda$ be a real number. If $b(E,X,\alpha) \ge \lambda$, then $b(F,X,\alpha) \ge \lambda .$
\end{lemma}

\begin{pf}
Let 
$$
\Delta_Y = f_*^{-1}\Delta_{X,\alpha} - \sum_{E'} a(E',X,\Delta_{X,\alpha})E'
$$ 
where the sum runs through all the $f-$exceptional divisors $E'$.  Then $f_*\Delta_Y = \Delta_{X,\alpha}$ and $K_Y+\Delta_Y \equiv f^*(K_X+\Delta_{X,\alpha}).$  Therefore, we can apply Lemma 2.30 of \cite{KM}, which gives,
$$
a(F,X,\Delta_{X,\alpha})=a(F,Y,\Delta_Y).
$$
Now, by Lemma \ref{discrep1}, we get,
$$
a(F,Y,\Delta_Y)=c-1-\left[ \sum a_i \cdot mult_{E_0}(f_*^{-1}D_i)-\sum_{E'}a(E',X,\Delta_{X,\alpha})\cdot mult_{E_0}E' \right]
$$
where $c=\mbox{codim}_YE_0$ and $\Delta_{X,\alpha}=\sum_i a_iD_i$.
Now,
$$
\Delta_{Y,\alpha}=\sum_i a_i(f_*^{-1}D_i)+\sum_{E'}\left(1-\frac{1}{e_{E'}}\right)E'
$$
where $e_{E'}$ is the degree of the cover on $E'$ induced by $\alpha \in Br(k(X)) \cong Br(k(Y))$.
Again, by Lemma \ref{discrep1},
$$
a(F,Y,\Delta_{Y,\alpha})=c-1-\left[\sum a_i \cdot mult_{E_0}(f_*^{-1}D_i) +\sum_{E'}\left(1-\frac{1}{e_{E'}}\right)\cdot mult_{E_0}E' \right].
$$
Thus, we get,
\begin{eqnarray*}
a(F,Y,\Delta_Y)-a(F,Y,\Delta_{Y,\alpha})&=&\sum_{E'}\left[a(E',X,\Delta_{X,\alpha})+\left(1-\frac{1}{e_{E'}}\right) \right]mult_{E_0}E' \\
&=&\sum_{E'} b(E',X,\alpha)mult_{E_0}E'.
\end{eqnarray*}
Since $b(E',X,\alpha) \ge 0,$ $b(E,X,\alpha) \ge \lambda$ and $mult_{E_0}E=1,$ we have,
$$
a(F,Y,\Delta_Y)-a(F,Y,\Delta_{Y,\alpha}) \ge \lambda .
$$
Since $a(F,Y,\Delta_{Y,\alpha}) \ge 0$ by hypothesis, we get $a(F,Y,\Delta_Y) \ge \lambda .$  But we proved earlier that $a(F,X,\Delta_{X,\alpha})=a(F,Y,\Delta_Y)$.  Thus, $a(F,X,\Delta_{X,\alpha}) \ge \lambda .$  This gives $b(F,X,\alpha)=a(F,X,\Delta_{X,\alpha})+\left(1-\frac{1}{e_F}\right) \ge \lambda .$
\end{pf}

\end{subsection}

\begin{theorem}
Any Brauer pair admits a terminal resolution.
\end{theorem}
\begin{pf}
Let $(X,\alpha)$ be a Brauer pair.   Using Hironaka's desingularization theorem \cite{Hiro}, we can assume that $X$ is smooth and $\Delta_{X,\alpha}$ is simple normal crossing, where $\Delta_{X,\alpha}$ is the boundary divisor on $X$ induced by $\alpha$.

If $(X,\alpha)$ is level 1 Brauer terminal then Lemma \ref{composition} with Lemma \ref{usualdiscrepancy} show that $(X,\alpha)$ is Brauer terminal.

If $(X,\alpha)$ is not level 1 Brauer terminal, then $X$ can be blown up using Lemma \ref{levelone} to obtain a level 1 Brauer terminal pair.  This pair is Brauer terminal, again by Lemmas \ref{usualdiscrepancy} and \ref{composition}.
\end{pf}

In summary, we have shown that given a Brauer pair $(X,\alpha),$ we can associate to it a Brauer pair $(Y,\alpha)$ with the following properties:
\begin{enumerate}
\item
$Y$ is nonsingular;
\item
There is a birational morphism $f:Y \longrightarrow X$;
\item
The boundary divisor $\Delta_{Y,\alpha}$ induced on $Y$ by $\alpha$ is simple normal crossing;
\item
The Brauer discrepancy of any exceptional divisor over $Y$ is positive.
 
\end{enumerate}

In short, any  Brauer pair $(X,\alpha)$ admits a terminal resolution $(Y,\alpha) \longrightarrow (X,\alpha)$.

\end{section}

\vskip 1cm

\noindent
{\bf Remark}.  Note that in our analysis, we restricted $\alpha$ to be in the 2-torsion part of the Brauer group $Br \ k(X)$.  If we require $\alpha$ to be in $Br_3 \ (k(X))$ instead, the analogous statement to the main result we proved here may not be true.  For example, consider a 3-fold $X$ with $\alpha \in Br_3\ k(X)$ that induces a simple normal crossing divisor $\Delta_{X,\alpha}$ that has the local form,
$$
\frac{1}{3} V(x_1) + \frac{1}{3} V(x_2) + \frac{1}{3} V(x_3).
$$
Suppose the 3-sheeted covers on $V(x_1)$ and $V(x_2)$ ramify on $V(x_1, x_2)$, but there is no ramification on $V(x_2,x_3)$ and $V(x_1,x_3)$.  This pair $(X,\alpha)$ cannot be improved by blowing up.   Now, let $E$ and $Y$ be the exceptional divisor and the variety generated respectively,  when $X$ is blown up along $V(x_1, x_3)$.   Then,
$$
b(E,X,\alpha) =2-\frac{1}{3} -\left(1-\frac{1}{3}+1-\frac{1}{3} \right) = \frac{1}{3},
$$
by Lemma \ref{bdiscrep1}.  Now blow up $Y$ along $E \cap V(x_1)$, and let $F$ be the exceptional divisor generated.  Then,
$$
b(F,Y, \alpha)=2-\frac{1}{e_F} - \left(1-\frac{1}{3} +1 -\frac{1}{3} \right)=\frac{2}{3} -\frac{1}{e_F}
$$
where $e_F$ is the degree of the cover on $F$ induced by $\alpha.$   We can show, using Brauer versions of Lemmas 2.29 and 2.30 of \cite{KM}, that 
$$
b(F,X, \alpha) =b(F,Y,\alpha)+b(E,X,\alpha).
$$
Then, we get
$$
b(F,X,\alpha)=\frac{2}{3}-\frac{1}{e_F} +\frac{1}{3}= 1-\frac{1}{e_F}.
$$
Note that it is not posible to determine the degree $e_F$ of the cover on $F$ induced by $\alpha$ without carrying out a detailed ramification computation involving roots of unity.  If it happens that $e_F=1$, then we get $B(F,X,\alpha)=0,$ indicating that the pair $(X,\alpha)$  may not admit a terminal resolution.  However, to determine the degree $e_F$ of the cover definitely, one must carry out ramification computations involving roots of unity.

\vskip 5mm

\noindent
{\bf Acknowledgement}:  This paper is based on the author's PhD thesis.  The author wishes to thank his supervisor, Colin Ingalls, for suggesting this problem and also for teaching the necessary background material.















\end{document}